\documentclass[11pt,reqno]{amsart}      
\usepackage{amssymb,amsthm}             
\usepackage{a4}                         
\numberwithin{equation}{section}        
\newcommand{\mnote}[1]{}  		

\renewcommand{\Re}{{\mathbb R}}         
\newcommand{\half}{\frac{1}{2}}         
\newcommand{\eps}{\epsilon}
\newcommand{\Ham}{\mathcal H}		
\newcommand{\tr}{\text{\rm tr}}		
\newcommand{\semidir}{\ltimes}          
\newcommand{\dev}{\mathcal D}		
\newcommand{\FF}{\mathcal F}		
		
\newcommand{\Vol}{\text{\rm Vol}}	

\newcommand{\ISO}{\text{\rm ISO}}

\newcommand{\SO}{\text{\rm SO}}

\newcommand{\iso}{\mathfrak{iso}}
\newcommand{\so}{\mathfrak{so}}

\newcommand{\Teich}{\text{\rm Teich}}	
\newcommand{\Tree}{\mathcal T}		
\newcommand{\Mink}{\mathbb M}		
\newcommand{\Hor}{\mathbf H}
\newcommand{\LL}{\mathcal L}		
\newcommand{\HH}{\mathbb H}

\newcommand{\CC}{\mathcal C}
\newcommand{\cover}{\tilde}
\newcommand{\univ}{\text{U}}
\newcommand{\tX}{\tilde X}

\newcommand{\genus}{\text{\rm genus}}		
\newcommand{\Energy}{\mathcal E}		

\newcommand{\conv}{\text{\rm conv}}

\theoremstyle{plain}
\newtheorem{thm}{Theorem}[section]
\newtheorem{conj}{Conjecture}
\newtheorem{claim}{Claim}

\newtheorem{lemma}[thm]{Lemma}

\title[CMC foliations of simplicial flat
spacetimes]{Constant mean curvature foliations of simplicial flat spacetimes}
\author[L. Andersson]{Lars Andersson}
\thanks{Supported in part by the Swedish Research Council, 
contract no.  R-RA 4873-307, the NSF,
contract no. DMS 0104402, and the Erwin Schr\"odinger Institute, Vienna.}
\address{Department of Mathematics\\
University of Miami\\
Coral Gables, FL 33124\\
USA}
\email{larsa\char'100math.miami.edu}

\date{July 25, 2003}

\allowdisplaybreaks[2]

\begin{document}
\begin{abstract}
Benedetti and Guadagnini \cite{benedetti:guadagnini:cosm} 
have conjectured that the constant mean curvature
foliation $M_\tau$ 
in a $2+1$ dimensional flat spacetime $V$ with compact hyperbolic Cauchy surfaces satisfies 
$\lim_{\tau \to -\infty} \ell_{M_\tau} = s_{\Tree}$, where $\ell_{M_\tau}$
and $s_{\Tree}$ denote the marked length spectrum of $M_\tau$ and the marked
measure spectrum of the $\Re$-tree $\Tree$, dual to the measured foliation
corresponding to the translational part of the holonomy of $V$,
respectively. We prove that this is the case for $n+1$ dimensional, $n \geq
2$, {\em simplicial} flat
spacetimes with compact hyperbolic Cauchy surface. A simplicial spacetime is 
obtained from the Lorentz cone over a hyperbolic manifold by deformations
corresponding to a simple measured foliation. 
 
\end{abstract}
\maketitle


\section{Introduction} 
In this paper we will consider maximal, globally hyperbolic, flat (MGHF) 
spacetimes $V$ of dimension $n+1$, $n \geq 2$, with compact Cauchy surface
$M$ of hyperbolic type, i.e. 
which admits a metric $g$ of constant sectional curvature $-1$. The main
result of the present paper implies, in the $2+1$ dimensional case, the proof of a conjecture of Benedetti and
Guadagnini \cite[Conj. 5.1]{benedetti:guadagnini:cosm}, 
see conjecture \ref{conj:benedetti} below,  in the special case of
simplicial flat spacetimes. A {\bf simplicial} 
spacetime is a flat spacetime which
can be obtained from the Lorentz cone over $M$, with metric $- d\rho^2 + \rho^2 g_0$
over $M$, by performing certain deformations relating to a weighted, finite
collection of nonintersecting compact simple totally geodesic hypersurfaces
$\LL = \{ (\Sigma_k, \ell_k), k = 1,\dots,m\}$, with weigths 
$\ell_k \in \Re$,  in $(M,g_0)$. In particular, a simplicial spacetime has a
compact Cauchy surface of hyperbolic type. Let $I_\ell = [0,\ell]$. 
The deformation corresponding to a single
such hypersurface $(\Sigma, \ell)$ corresponds to gluing in a
``wedge'' spacetime $W_\ell = 
W \times I_\ell$ in place of the Lorentz cone $W$ over
$\Sigma$. In case $n=2$, $\LL$ is a ``weighted
multicurve'', or simple measured foliation with compact simple geodesic leaves 
(by \cite{zeghib:hyp:geod}, every totally geodesic
hypersurface of a compact hyperbolic manifold of dimension $n \geq 3$ is
compact).  
Simple measured foliations with compact leaves are dense in the space of
all measured foliations. 

A MGHF spacetime $V$ with compact Cauchy
surface $M$ of hyperbolic type 
is globally foliated by CMC hypersurfaces $M_\tau$ with $\tau$ taking all
values in $(-\infty,0)$, see \cite{andersson:flatcmc}. 
Further, the scale free version $\frac{\tau^2}{n^2} g_\tau$
of the induced metric $g_\tau$ on $M_\tau$ converges in the expanding direction,
as $\tau \nearrow 0$, to a metric of constant sectional curvature $-1$. In
case $n \geq 3$, this metric is the unique hyperbolic metric on $M$, while in
case $n=2$, this metric corresponds to a point in the Teichmuller space
$\Teich(M)$ of $M$. This is a partial generalization of the results for the
case $n=2$ proved in 
\cite{andersson:etal:2+1grav}. In that paper it was also proved that in the
direction $\tau \searrow -\infty$, towards the singularity,  the Teichmuller
class of the induced metric on $M_\tau$ diverges, in the sense that it 
leaves every compact subset of
$\Teich(M)$, as $\tau \searrow -\infty$. 
This is proved by showing that the Dirichlet energy $\Energy$, 
which is a
proper function on $\Teich(M)$ \cite[\S 3]{schoen:yau:incompressible}, see
also \cite{tromba:book}, 
diverges as $\tau \searrow -\infty$. However, the work in
\cite{andersson:etal:2+1grav} does not give a detailed picture of the
geometry of the CMC hypersurfaces $M_\tau$ for $\tau \searrow -\infty$. It is
the purpose of this paper to study the detailed asymptotic 
behavior of the geometry of
$M_\tau$ in the case when $V$ is simplicial. 

\subsection{Flat spacetimes, earthquakes and $\Re$--trees}
A time oriented MGHF spacetime $V$ with oriented Cauchy surface 
$M$ may be viewed as an $\ISO^+(n,1)$
geometric structure, and as such is described by the holonomy
representation $\alpha$ 
of the fundamental group $\pi_1(M)$ in $\ISO^+(n,1)$. The
decomposition $\ISO(n,1) = \SO(n,1)\semidir \Re^{n+1}$ leads to a
decomposition $\alpha(\gamma) = Q(\gamma) + t_\gamma$ where $Q$ is the linear
part of the holonomy and the translational part $t_\gamma$ is a
cocycle with values in $\Re^{n+1}$. The linear part $Q$ of $\alpha$ corresponds to a hyperbolic
structure on $M$, i.e. a point in $\Teich(M)$. 
Let $\Gamma = \alpha(\pi_1(M)) \subset
\ISO^+(n,1)$. The moduli space of MGHF spacetimes with
Cauchy surface $M$ is homeomorphic to an open ball 
in the Zariski tangent space 
$H^1(\Gamma, \iso^+(n,1)_{\text{\rm Ad}})$, 
see \cite{mess:const:curv,andersson:flatcmc}. We denote by $\iso$ and $\so$
the Lie algebras of $\ISO$ and $\SO$, respectively, and $\text{\rm Ad}$
indicates that $\Gamma$ acts by the adjoint representation. 

In case $n=2$, the
dimension of this space is $12 \genus(M) -12$, twice the dimension of
Teichmuller space, while for $n \geq 3$, by
Mostow rigidity we have $H^1(\Gamma, \iso^+(n,1)_{\text{\rm Ad}}) = H^1(\Gamma,
\Re^{n+1}_{\text{vec}})$. For $n=3$, the
dimension of the moduli space of flat spacetimes is the same as that of the
space of flat conformal structures, 
$H^1(\Gamma,\Re^{3+1}_{\text{vec}}) = H^1(\Gamma, \so(4,1)_{\text{\rm Ad}})$ 
\cite[\S 11]{kapovich:millson:deform}. 
In case $n=2$, the translational part
$t$ corresponds to a measured foliation $\FF$ of $M$. Given a measured
foliation $\FF$ of a hyperbolic surface $M$, there is a unique geometric real
tree ($\Re$-tree) $(\Tree,d)$ dual to it. The fundamental group $\pi_1(M)$
acts in a natural way on $\Tree$. 

Returning to the case of general dimension $n$, 
the development $\dev(\univ(V))$ of the universal cover $\univ( V)$ in
$\Mink^{n+1}$, the $n+1$ dimensional Minkowski space, 
is a convex open subset, the boundary $\Hor$ 
of which is the Cauchy
horizon of $\univ(V)$. From general results in causality theory, 
$\Hor$ is a weakly spacelike $C^0$ hypersurface. The fundamental group
$\pi_1(M)$ acts on $\dev(\univ(V))$ and the action extends to $\Hor$. 
The Lorentz structure of
$\Mink^{n+1}$ induces a degenerate distance function 
$d$ on $\Hor$. In the simplicial case, if we identify points of $\Hor$
under the equivalence relation $\sim$ defined by $p \sim q$ if $d(p,q) = 0$,
the action of $\pi_1(M)$ drops to the quotient 
$\Hor/ \!\! \sim$.
The
metric space $\Hor/ \!\! \sim$ can be identified with $(\Tree,d)$, 
which in this case is simplicial. 
It would be interesting to know whether the above direct construction of
$\Tree$ works for general MGHF spacetimes. 

The cosmological time
function $\rho(p)$, see  \cite{andersson:etal:cosmtime},
on $\dev(\univ(V)) \subset \Mink^{n+1}$, 
defined as the maximal Lorentz length of
any past directed causal curve starting at $p$ in $\dev(\univ(V))$ is
regular, i.e. it is everywhere finite and $\rho \to 0$ along every past
directed inextendible causal curve. This construction drops to the quotient
$V$. 
Starting from the 
work of Mess \cite{mess:const:curv}, 
Benedetti and Guadagnini \cite{benedetti:guadagnini:cosm}
showed that in case $n=2$, 
the induced geometry of the level sets of the cosmological time
function $\rho$ introduced in \cite{andersson:etal:cosmtime} realize the 
Thurston earthquake deformation, in the sense that the curve in Teichmuller
space defined by the Teichmuller class of
the induced geometry of
the level sets $M_\rho$ of
the cosmological time function corresponds to the Thurston earthquake flow,
defined with respect to the hyperbolic structure given by $Q$ and the 
measured foliation $\FF$, 
see 
\cite[Prop. 4.27 and \S 4.6]{benedetti:guadagnini:cosm}. In particular, as
$\rho \to \infty$, the Teichmuller class of $M_\rho$ converges to the
hyperbolic surface $M$ with holonomy $Q$, while as $\rho \to 0$, the
geometry of the universal cover $\univ(M_\rho)$ converges in the Gromov sense
to an $\Re$--tree $\Tree$, 
determined by the translational part $t$ of the holonomy $\alpha$ of $V$. 

The $\Re$--tree $\Tree$ can be identified with a point on the 
Thurston boundary of Teichmuller space. To explain this fact 
this we need the notions
of {\bf marked length spectrum} and {\bf marked measure spectrum}, which we
introduce following 
\cite[\S 4.5]{benedetti:guadagnini:cosm}. Let 
$(\tX, d)$ be a metric space with an action $\alpha$ of $\pi_1(M)$ and $X =
\tX/\pi_1(M)$. Let $\CC$ be the space of conjugation classes of 
$\pi_1(M) \setminus \{1\}$. Then $\CC$ can
be identified with the space of nontrivial homotopy classes of closed curves
on $M$. For $c \in \CC$, the marked length spectrum $s_X(c)$ is defined as 
$s_X(c) = \inf_{p \in \tX} d(p,\alpha(p))$. In case $X$ is homeomorphic to
$M$, this corresponds to the shortest length of closed curves in $c$, and is
denoted by $\ell_X$. In particular, by letting $X$ vary among the
hyperbolic structures on $M$, $s_X$ gives a map $s_X : \Teich(M) \to
\Re_{\geq 0}^{\CC}$, which is strictly positive. 

On the other hand, 
in case $\tX = \Tree$, $s_{\Tree}(c)$ can be expressed in terms of the measured
foliation $\FF$ of $M$ dual to $\Tree$ as the minimal transverse measure
realized by the curves in $c$. This gives the marked measure spectrum
$I_{\FF}$. If $\FF$ is a simple measured foliation $\LL$ 
with
compact leaves, then $I_{\LL}(c)$ is defined in terms of 
the geometric intersection number of
the curves in $c$ with $\LL$. 

This extends the notion of length spectrum to the degenerate case. In the 
2 dimensional case, the 
image of $\Teich(M)$ under $\ell_X$ is
homeomorphic to the open ball in 
$\Re^{6\genus(M)-6}$. The boundary consists of 
degenerate geometries corresponding to projective rays in the image of the
space of measured foliations under $I$. This is the Thurston boundary of
Teichmuller space. 
The convergence of marked spectra can be understood as 
convergence of metric spaces in the Gromov sense, see \cite[Remark 4.24,
point 3)]{benedetti:guadagnini:cosm}. 

By \cite{benedetti:guadagnini:cosm}, 
the foliation $M_\rho$ gives an analytic curve in Teichmuller space
connecting the interior point $(M,Q)$ to the point on the Thurston boundary
corresponding to $\Tree$. Thus,
the spacetime geometry allows us to recover all the information about the
holonomy in a concrete way. 
In the particular case of a $2+1$ dimensional 
flat simplicial spacetime, defined by a
hyperbolic surface $M$ and a simple measured lamination with compact leaves
$\LL$ on $M$, the Teichmuller
class of the level sets $M_\rho$ of the cosmological time function sweep out
a curve corresponding to the Fenchel--Nielsen deformation of $M$
obtained by twisting $M$ along the closed geodesics $\Sigma_k$ of $\LL$, and
the geometry on $\univ(M_\rho)$ converges in the Gromov sense to the
simplicial tree $\Tree$ dual to $\LL$. 
We refer to
\cite{mess:const:curv,benedetti:guadagnini:cosm,otal:hyp} for background on the concepts
discussed above. 

The conjecture of Benedetti and Guadagnini
can now be stated as follows: 

\begin{conj}[\protect{\cite[Conj. 5.1]{benedetti:guadagnini:cosm}}] 
\label{conj:benedetti} 
Let $V$ be a $2+1$ dimensional MGHF spacetime
with compact Cauchy surface of genus $\geq 2$, and let $M_\tau$ be the foliation of
$V$ by constant mean curvature hypersurfaces with mean curvature $\tau$. Then 
\begin{enumerate}
\item \label{point:infty} $\lim_{\tau \to -\infty} \ell_{M_\tau} = s_{\Tree}$,
\item \label{point:0} $\lim_{\tau \to 0} \ell_{M_{\tau}}/\tau = \ell_M$.
\end{enumerate}
\end{conj}
Point (\ref{point:0}), which states that the scale--free geometry on $M_\tau$
converges to the hyperbolic geometry $(M,g)$ corresponding to the holonomy
$Q$ in the expanding direction $\tau \to 0$, 
follows for $n \geq 2$ from the work in \cite{andersson:flatcmc}. 
In this paper we will prove that the statement corresponding to point
(\ref{point:infty})  is true for
{\em simplical} 
MGHF spacetimes with compact Cauchy surface of hyperbolic type, of
general dimension $n \geq 2$. We can state the main result of this paper as
follows, see Theorem \ref{thm:simpconv}.  
\begin{thm}\label{thm:main}
Let $V$ be an $n+1$ dimensional simplicial spacetime and let $M_\tau$ be
the foliation of $V$ by constant mean curvature hypersurfaces with mean
curvature $\tau$. Then $lim_{\tau \to -\infty} \ell_{M_\tau} = s_{\Tree}$. 
\end{thm}

Recall that in case $n =2$, the simple measured foliations with compact
leaves are dense in the space of all measured foliations. It is therefore
natural to conjecture that the 
result proved here will yield the general case by a limit argument. We will
not consider this problem here. 

Our results here hold for simplicial flat spacetimes 
in general dimension $n+1$, $n \geq 2$. The relation of the
case of simplicial flat spacetimes to the general case can be 
expected to be quite complicated in
higher dimensions. In fact, Scannell \cite{scannell:deform} showed
there are nonrigid compact hyperbolic 3--manifolds (i.e. ones with
$H^1(\Gamma, \so(4,1)_{\text{\rm Ad}}) \ne \{0\}$) which have no immersed
totally geodesic hypersurfaces. Therefore, it is not clear if the 
$2+1$ dimensional picture described above generalizes to the higher
dimensional case. It is an interesting open problem to describe the
asymptotics of both 
the foliation by level sets of the cosmological time function and of the 
CMC foliation of general higher dimensional flat spacetimes. 

\mnote{
{\bf Holonomy of $V_{\ell}$:}
cf.  Lafontaine, he computed f so that the Codazzi tensor 
$\tilde p$ is given by 
$\tilde p = \nabla^2 f - f g$ in terms
of the distance to the hyperplane corresponding to $\Sigma$. 
One would like to have 
a nice formula for the cocycle $\sigma$ corresponding 
to the wedge over $\Sigma$}

One of the main ideas in the work of Benedetti and Guadagnini is that the
foliation by level sets of the cosmological time function realizes in a
natural way the earthquake deformation of Thurston, with respect to the
measured foliation defined by the translational part of the holonomy of the
spacetime. It is an interesting problem to understand the corresponding
picture in the higher dimensional case. As discussed below, the level sets of
the cosmological time function in a flat simplicial spacetime have conformally flat induced geometry.
Recall the coincidence 
$H^1(\Gamma , \Re^{3+1}_{\text{\rm vec}}) = H^1(\Gamma,\so(4,1)_{\text{\rm Ad}})$,
which holds in dimension 3 only.
This raises the 
possibility that the cosmological time foliation in a general $3+1$
dimensional flat spacetime gives a parametrization of the deformation space
of flat conformal structures on $M$, in a way analogous to the
$2+1$ dimensional case described above. 
\mnote{It would be interesting to know whether the family $M_\rho$ in a given
$V$ gives a ``bending'' deformation of the hyperbolic structure $M$ in the
space of flat conformal structures.}
Since the moduli space of MGHF  $n+1$ dimensional spacetimes with compact
hyperbolic Cauchy surface is a manifold
\cite{andersson:flatcmc}, if this relation is true, it would imply
the conjecture of Kapovich 
\cite{kapovich:deform}, 
that the space of flat conformal structures on a compact hyperbolic
manifold is smooth in dimension 3. 

\section{CMC hypersurfaces in wedge spacetimes}\label{sec:wedge}
Let $(M,g)$ be compact a compact hyperbolic manifold of sectional curvature
$-1$, 
with compact totally geodesic embedded hypersurface
$\Sigma$, and denote the induced hyperbolic (if $n \geq 3$) 
metric on $\Sigma$ by $h$. 
Let $V = (0,\infty)\times M$ be the flat
Lorentz cone over $M$ with metric 
$$
ds^2 = - d\rho^2 + \rho^2 g, \qquad \rho \in (0,\infty) .
$$
For $\ell > 0$, the {\bf wedge spacetime} $V_\ell$ is 
$V$, with the cone over $\Sigma$ replaced by the wedge
$W_{\ell}$ of width $\ell$, given by 
$$
W_{\ell} = (0,\infty) \times \Sigma \times I_\ell ,
$$
with metric 
$$
- d\rho^2 + \rho^2 h + dr^2 , \qquad (\rho,r)\in (0,\infty)\times
  I_\ell  .
$$
$V_\ell$ is a MGHF simplicial 
spacetime which is a deformation of $V$. The above type of deformation was
called {\bf elementary} in \cite{benedetti:guadagnini:cosm}. 
It will be useful to pass to the covering of these spacetimes defined
w.r.t. the fundamental group $\pi_1(\Sigma)$. We use notation of the form 
$\cover V_\ell$ or $\univ_\Sigma(V_\ell)$ for this cover, while $U(V)$
denotes the universal cover. Let $I^{n+1}_+(\{0\})$ denote the interior of the
future light cone of the origin in $\Mink^{n+1}$. Then 
$\cover W_\ell = I^n_+(\{0\}) \times I_\ell$.
In coordinates $t,y,r$, $\cover W_\ell$ is the set
$-t^2 + |y|^2 < 0$, with metric 
$$
-dt^2 + dy^2 + dr^2 .
$$
The level sets $\cover M_{\rho}$
of $\rho$ in $\cover V_{\ell}$ has metric $\rho^2 g$ in $\cover ( V \setminus
\Sigma) $ and
metric $\rho^2 h$ in $\univ (\Sigma \times I_\ell) $. 
This metric is $C^1$ but not $C^2$, the second derivatives being bounded but
not continuous, and it is conformally flat. 
To see this explicitely,
note that in the Gauss foliation based on $\Sigma$, the
metric $g$ can be written in the form the form 
\begin{equation}\label{eq:gaussconf}
g = \cos^{-2}(v)(dv^2 + h), \qquad v \geq 0,
\end{equation}
where $v = 0$ at $\Sigma$. 
In case $n=2$, the wedge metric is flat, and the above form of $g$ shows that
it is conformally flat. 
Next we consider the 3 dimensional case. The Cotton tensor is 
$
C_{ijk} = 2 \nabla_{[k} ( R_{j]i} - \frac{1}{4} R g_{j]i} ) .
$
The vanishing of $C_{ijk}$ characterizes local conformal flatness in
dimension 3. 
In case
$n=3$, the wedge metric for $\rho = 1$, 
$g = h + dr^2$ has
$C_{ijk} = 0$, so $g$ is conformally flat. From equation (\ref{eq:gaussconf})
we see that off the wedge, $g$ is conformal to a metric of the same form as
the wedge metric and hence is conformally flat. 
Finally, in case $n > 3$, the metric $h + dr^2$ has 
nonvanishing Weyl tensor so it is not conformally flat. 
\bigskip

\subsection{Mean curvature of $M_{\rho}$}
The second fundamental form is 
$
K = - \half \partial_\rho g(\rho) .
$
On $M \setminus \Sigma$ we have $K = - \rho^{-1} g$, while on $\Sigma \times 
I_\ell$
we have $K = - \rho^{-1} h \oplus 0$. 
The mean curvature $\tau = \tr K$ is given by $\tau = - n/\rho$ on
$M\setminus \Sigma$ while on $\Sigma\times I_\ell$, $\tau =
-(n-1)/\rho$. This means in particular that if we choose $\rho_0, \rho_1$ so
that 
$$
- \frac{n-1}{\rho_0} < \frac{-n}{\rho_1} ,
$$
then 
$$
\max \left ( \tau \bigg{|}_{M_{\rho_0,\ell}} \right ) < 
\min \left ( \tau \bigg{|}_{M_{\rho_1,\ell}} \right ) .
$$
This shows that the level sets $M_{\rho}$ are barriers, in the sense of 
\cite{andersson:etal:maxprin}, for the mean
curvature equation in
$V_{\ell}$, which using the argument of Gerhardt \cite{gerhardt:CMC} 
gives an easy proof that the wedge space--times $V_{\ell}$ 
are globally foliated by CMC hypersurfaces. 
The function $\rho$ defined above is the cosmological time
\cite{andersson:etal:cosmtime} of $V_{\ell}$.

\subsection{CMC hypersurfaces}
Now consider the CMC hypersurfaces 
$M_{\tau}$ of mean curvature $\tau < 0$, 
in the unique global CMC foliation of $V_{\ell}$. 
We will scale  
$V_{\ell}$ by a factor $\lambda^2$, the rescaled metric is 
$g' = \lambda^2 g$. 
This has the effect of scaling $\tau$ to $\lambda^{-1} \tau$. 
We shall choose 
$$
\lambda = |\tau|/(n-1),
$$ 
so that the rescaled version of the
hypersurface $M_{\tau}$ 
has mean curvature $- (n-1)$, and  
consider the limit as $\tau \to - \infty$, i.e. as 
$\lambda \to \infty$. 

The scaling changes $V_\ell$ to $V_{\lambda\ell}$, in particular the wedge in
$V_{\lambda\ell}$ is $W_{\lambda\ell}$, which after a change of coordinates
$\rho' = \lambda \rho$, $r' = \lambda r$, has metric of the form 
\begin{equation}\label{eq:prim}
- (d\rho')^2 + {\rho'}^2 h_{ij}dx^i dx^j 
+ {dr'}^2 , \qquad r' \in   I_\ell  .
\end{equation}
where $x^i$, $i=1,\dots,n-1$ is a coordinate system on $\Sigma$. On $\tilde
W_{\lambda\ell}$ we also have the 
scaled Minkowski coordinate system $(t',y',r')=\lambda (t,y,r)$, with metric 
$$
-(dt')^2 + (dy')^2 + (dr')^2 ,
$$
so that ${\rho'}^2 = {t'}^2 - |y'|^2$. 
We see from this that the scaling has the effect of stretching the
wedge $W_{\ell}$ to the wedge $W_{\lambda\ell}$ of width $\lambda \ell$. 
We denote the unique CMC hypersurface in $V_{\lambda\ell}$ with mean
curvature $-(n-1)$ by
$M_{\lambda}$.
Let $u_{\tau}$ and $u_{\lambda}$ denote the height functions of $M_\tau$ and
$M_\lambda$ with respect to the time function $\rho$, 
defined by $u_\tau = \rho \big{|}_{M_\tau}$ and 
$u_\lambda = \rho' \big{|}_{M_\lambda}$ and let $\cover u_\tau, \cover
u_\lambda$ denote the corresponding lifts. Similarly, let $v_\tau = t
\big{|}_{\cover M_\tau}$ and $v_\lambda = t' \big{|}_{\cover M_\lambda}$. 

In view of the mean curvature of the level sets of $\rho$, we have by the
maximum principle,  
$\lambda^{-1} \leq u_\tau \leq  \lambda^{-1}n/(n-1)$, and 
$1 \leq u_\lambda \leq n/(n-1)$. 
The mean curvature of $M_\lambda$ is $-(n-1)$, and hence the 
derivative bounds for constant mean curvature hypersurfaces
\cite{cheng:yau:maximal,treibergs:cmc} apply to $v_\lambda$. It 
follows that there is a subsequence of $u_{\lambda}$ which converges
uniformly in $C^3$ on compacts to a hypersurface $M_{\infty}$ with mean
curvature $-(n-1)$ 
in $W_{\infty}$ where $W_{\infty}$ is the
Kasner type space--time 
$(0,\infty) \times \Sigma \times \Re$ with metric 
$$
- d\rho^2 + \rho^2 h + dr^2, \qquad - \infty < r < \infty
$$
This space--time is the product of the flat Lorentz cone over $\Sigma$ with a
line. 

The conclusion so far is that the limiting hypersurface is an entire CMC
hypersurface in $W_{\infty}$, with mean curvature $-(n-1)$. Further, 
due to the fact that $1 \leq u_\lambda \leq n/(n-1)$, 
$M_{\infty}$ lies between the barriers $\rho = 1$ and $\rho =
n/(n-1)$. 
In fact, as well will now prove, a surface $M_{\infty}$ with these properties
splits as a product. We state this as the following 
\begin{claim}
Let $M$ be an entire CMC hypersurface of mean curvature $-(n-1)$ in
$W_{\infty}$, bounded from above and below by by the barriers $N_1, N_2$
\begin{align*}
N_1 &= \{\rho = \rho_1\}, \rho_1 \leq 1 \\
N_2 &= \{\rho = \rho_2\}, \rho_2 \geq  n/(n-1) 
\end{align*}
Then $M$ splits as $M=\Sigma\times \Re$ and $M$ coincides with the level set
$\rho =1$. 
\end{claim}
We will prove the claim as a special case of a more general splitting
theorem.
\begin{thm} \label{thm:split}
Let $W = (0,\infty) \times \Sigma^n \times \Re^k$, with metric 
$$
ds^2 = - d\rho^2 + \rho^2 h + (dz^1)^2 + \cdots (dz^k)^2 .
$$
Let $M$ be an entire CMC hypersurface in $W$ with mean curvature $-n$,
bounded between the barrier surfaces 
\begin{align*}
N_1 &= \{\rho = \rho_1\}, \quad \rho_1 \leq 1 , \\
N_2 &= \{\rho = \rho_2\}, \quad \rho_2 \geq n/(n-1) . 
\end{align*}
Then $M$ splits as a product $M = \Sigma \times \Re^k$ with metric 
$$
h + (dz^1)^2 + \cdots (dz^k)^2 ,
$$
and $M$ coincides with the level set $\rho =1$. 
\end{thm}
\begin{proof}
Recall that the universal cover of the Lorentz cone over $\Sigma$ is
$I^{n+1}_+(\{0\})$, the interior of the future light cone in the
$n+1$ dimensional Minkowski space $\Mink^{n+1}$.
The universal cover $p: \tilde \Sigma \to \Sigma$ induces the universal cover
$p: \tilde W \to W$, with $\tilde W = I^{n+1}_+(\{0\}) \times \Re^{k} \subset
\Mink^{n+1+k}$, the $n+1+k$ dimensional Minkowski space.  
Let $\tilde M$ be the lift of $M$ to $\tilde W$. Then
$\tilde M$ is an entire CMC hypersurface in $\tilde W$ which we therefore may
think of as a CMC hypersurface in $\Mink^{n+1+k}$, the $n+1+k$ dimensional
Minkowski space.
Introduce coordinates $(t,y^1,\dots,y^n,z^1,\dots,z^k)$ on $\Mink^{n+1+k}$. We
will use the notation $x = (y,z)$. Let
$|y|^2 = (y^1)^2 + \cdots + (y^n)^2$, $|z|^2 = (z^1)^2 + \cdots + (z^k)^2$,
and define the function $\tilde \rho$ on $\Mink^{n+1+k}$ by 
$$
\tilde \rho^2 = t^2 - |y|^2 .
$$
Then $\tilde N_i = \{ \tilde \rho = \rho_i \}$, $i=1,2$ are the universal
covers of $N_1, N_2$, and from the maximum principle and the 
assumptions of the theorem it follows that
$\tilde M$ is bounded between $\tilde N_1 $ and $\tilde N_2$. 

We will now make use of some results of Choi and Treibergs
\cite{treibergs:choi}. 
The conclusion of \cite[\S 4]{treibergs:choi} can be summarized as follows. 
Let $v$ be the height function of a $\tau \ne 0$ CMC hypersurface $M \subset
\Mink^{n+k+1}$, $v = t \big{|}_{M}$.
Let $V_v$ be the positive homogenous of degree one function
defined by 
$$
V_v = \lim_{r \to \infty, \  r > 0} \frac{u(rx)}{r} .
$$
By \cite[Lemma 4.6]{treibergs:choi}, the tangent cone to $V_v$ at $0$,
$\chi_{V_v}$ is given by 
$$
\chi_{V_v}(0) = \conv(L_v),
$$
the convex hull of some closed subset $L_v$ in $\HH^{n+k}(\infty)$. Here
$\HH^{n+k}$ may be identified with the unit ball in $\Re^{n+k}$ with
coordinates $(y,z)$, so that $\HH^{n+k}(\infty)\cong S^{n+k-1}$. 
Let $E^n = \{ (y,z)\in \Re^{n+k}: z = 0\}$. 
By \cite[Lemma 4.3]{treibergs:choi}, cf. proof of \cite[Lemma 4.6]{treibergs:choi} 
$$
V_v = \sup_{\xi \in L_v} x \cdot \xi \, .
$$
We now make the following 
\begin{claim}
$L_v \subset S^{n+k-1} \cap E^n$. 
\end{claim}
If
this holds then the splitting theorem \cite[Theorem 4.8]{treibergs:choi} shows that in
fact $M$ splits as $M^n \times \Re^k$. 
The claim will follow if $V_v(0,z) =
0$ for all $x \in \Re^k$.   
Let $v_2$ be the height function of the future barrier 
$N_2$. 
By construction, $v_2(y,z) = w_2(y)$, in
particular $w_2$ is independent of $z$. It follows that $V_{v_2}(0,z) = 0$. 
Since $v(x) \leq v_2(x)$ we have $V_v(0,z) \leq V_{v_2}(0,z)$ and hence 
$V_v(0,z) = 0$, which proves the claim. 
It follows that $L_v \subset S^{n+k-1} \cap E^n$ and hence $M$ splits as a
metric product $M = M^n \times \Re^k$, where $M^n$ is a CMC hypersurface of
$n+1$ dimensional Minkowski space $\Mink^{n+1}$. 

Applying this result to the universal cover $\tilde M$ we see that the
splitting also applies to $M$ and the conclusion is that $M$ splits as 
$M = \Sigma \times
\Re^k$. By assumption, the mean curvature of $M$ is $-(n-1)$ which due to the
split of $M$ implies that $M = \{\rho = 1\}$. 
\end{proof}
Going back to the limiting process, we see that we have proved that
$M_{\lambda}$ converges on compacts to the metric product $\Sigma \times \Re$. 
In terms of the
height function $u_{\lambda}$ we have proved
\begin{lemma}\label{lem:conv}
$u_\lambda$ converges uniformly in $C^2$ 
on compacts in $W_{\lambda\ell}$ to the constant function $1$. 
\end{lemma}
By the barrier construction we have 
\mnote{here we should have strict inequality by the strong maximum principle} 
\begin{align}
1 &< u_{\lambda} < \frac{n}{n-1} , \label{eq:ulambdabound} \\
\intertext{or} 
\lambda^{-1} &< u_{\tau} < \lambda^{-1} \frac{n}{n-1} . \label{eq:utaubound}
\end{align}
The second fundamental form $K$ of
$M_{\tau}$ satisfies $K \leq 0$ with our conventions. 
This means that the height functions w.r.t. $t$ and $t'$, $v_\tau$ and
$v_\lambda$ are convex, cf. \cite[Prop. 1.1]{treibergs:choi}, 
in particular 
$$
\partial^2/\partial {r}^2 v_{\tau}(y,r) \geq 0 .
$$
We have 
\begin{equation}\label{eq:vtau-utau}
v_\tau^2  = \cover u_\tau^2 + |y|^2 .
\end{equation}
From the above, $\cover u_\tau$ varies by at most $\lambda^{-1}$ which means
that $|v_\tau(y,r_1) - v_\tau(y, r_2)| \leq \frac{1}{(n-1)\lambda}$.
Further, if we restrict to one fundamental domain of $\cover M_\tau$, the
projection on the $y$-variables is bounded by $C\lambda^{-1}$.
 
We shall need the following elementary calculus Lemma. 
\begin{lemma} Let $f: [a , b]$ be a convex $C^2$ 
function which takes
values in $[0,\Delta]$. Then for any $\eps > 0$, $\eps < (b-a)/2$.  
the estimate 
$$
|f'(x)| \leq \frac{\Delta}{\eps} ,
$$
holds in the interval $[a+\eps, b-\eps]$. 
\end{lemma}
Let $I_{\ell,\eps} = (\eps, \ell - \eps)$. 
Note that $u_\tau \partial_{r} u_\tau = v_\tau \partial_{r} v_\tau$, 
and hence in view of the above mentioned bound on $y$ in a fundamental domain
of $\cover M_\tau$, and the lower bound on $u_\tau$, 
the Lemma applies to 
applies to $\partial_r u_\tau$ to give an estimate of the form 
$$
|\partial_r u_\tau | \leq \frac{C}{\lambda \eps} , 
\qquad \text{ for } r \in I_{\ell,\eps} .
$$
The derivative bounds give $|D' v_\lambda| \leq C$,
$|{D'}^2 v_\lambda | \leq C $, $|{D'}^3 v_\lambda|\leq C$ over
compacts. Taking into account the boundedness of the fundamental domain of
$\pi_1(\Sigma)$ 
in $\cover\Sigma$ and consequently in $M_\lambda$, and the relation of
$u_\lambda$ to $v_\lambda$, we have the corresponding bounds for
$u_\lambda$. The same bounds hold also in terms of the coordinates $x,r'$, $r
= \lambda r'$ on
$\Sigma \times \lambda I_\ell$. 

Now we consider $u_\tau$, and note that this is just a rescaling of
$u_\lambda$ by a factor $1/\lambda$, 
$$
u_\tau (x,r') = \lambda^{-1} u_\lambda(x,r') .
$$
This gives, in view of the fact that the $x$--coordinate does not scale,  
\begin{equation}\label{eq:Dkutau}
|D_x^k u_\tau| \leq C/\lambda, \quad k=1,2,3 .
\end{equation} 
From Lemma \ref{lem:conv} we have
also $|D_r u_\tau| \leq \frac{C}{\lambda\eps}$, for 
$r \in I_{\ell,\eps}$. Without the use of 
the Lemma, we would just have an estimate of the form $|D_r u_\tau|\leq C$. 
\begin{lemma} \label{lem:Dxu0}
Fix $(x_0, r_0) \in \Sigma \times I_\ell$. Then 
$$
\lim_{\lambda \to \infty} \lambda D_x u_\tau (x_0, r_0) \to 0 .
$$
\end{lemma}
\begin{proof}
We have $\lambda D_x u_\tau (x,r) = D_x u_\lambda(x,r')$. 
Let $r'_0 = \lambda r_0$, $\bar r = r' - r'_0$ and $\bar u_\lambda(x,r') =
u_\lambda(x, r'-r'_0)$. This has the effect of translating $r'_0$ to $0$. 
The derivative bounds apply to $\bar u_\lambda$ and hence also the conclusion
of Theorem \ref{thm:split}, which implies that $\bar u_\lambda \to 1$ in
$C^2$ on compacts. 
The result follows. 
\end{proof}
To compute the induced metric on $M_{\tau}$ we work in coordinates  
$(x,r)$, $x = (x^1, \dots, x^{n-1})$ on $\Sigma \times I_\ell$, 
and define the map $\Phi_\tau : \Sigma \times I_\ell \to W_{\ell}$, by 
$$
\Phi_\tau(x,r) = (u_{\tau} (x,r) , x,r) .
$$
Then the image of $\Phi_\tau$ is precisely $M_\tau \cap
W_{\ell}$. 
Let the indices $i,j$ run over $1,\dots,n-1$ and let the index $n$ correspond
to the coordinate $r$. Pulling back the metric $-d\rho^2 + \rho^2 h + dr^2$
by $\Phi_\tau$ gives 
$$
g_\tau = u_\tau^2 h \otimes 1 - du_\tau \otimes du_\tau , \qquad \text{ in }
M_\tau \cap W_\ell ,
$$
which shows that $g_\tau \leq u_\tau^2 h \otimes 1$ as quadratic forms. 
From this follows
$$
\det g_\tau \leq u_\tau^{2(n-1)} \det h , \qquad \text{ in } M_\tau \cap
W_\ell .
$$
Similarly, we have 
$$
\det g_\tau \leq u_\tau^{2n} \det g, \qquad \text{ in } M_\tau \setminus
W_\ell ,
$$
where $g$ is the hyperbolic metric on $M$. From this follows in particular
that 
\begin{equation}\label{eq:Mtau-Well}
\lim_{\lambda \to \infty} \lambda^{n-1} \Vol(M_\tau \setminus W_\ell) = 0 .
\end{equation}
We have in view of the fact that $\lambda u_\tau \leq n/(n-1)$, 
\begin{equation}\label{eq:lameps}
\lambda^{n-1} \int_{\Sigma \times (I_\ell \setminus I_{\ell,\eps})}
\sqrt{\det g_\tau } dx dr
\leq C \eps \Vol(\Sigma) .
\end{equation}
First consider the case $n=2$. Then $\Sigma$ is $1$ dimensional with metric
$h dx^2$, 
and the explicit form of $\det g_\tau$ is
$$
\det g_\tau = [1-(\frac{\partial u_\tau}{\partial r})^2 ] u_\tau^2 h  -
(\frac{\partial u_\tau}{\partial x} )^2 \leq u_\tau^2 h .
$$
Here we may take $h\equiv 1$ by choosing $x$ to be the arclength parameter on
$\Sigma$. 
By Theorem \ref{thm:split}, $\lambda u_\tau \to 1$, and by Lemma
\ref{lem:Dxu0} $\lambda \partial u_\tau/\partial x \to 0$, pointwise as $\lambda \to
\infty$. The dominated convergence theorem now shows 
$$
\lambda \int_{\Sigma \times I_{\ell,\eps}} \sqrt{\det g_\tau} dx dr = (\ell
-2\eps) L(\Sigma) .
$$
where $L(\Sigma)$ denotes the length of $\Sigma$. Since $\eps > 0$ is
arbitrary, we conclude 
$$
\lim_{\lambda \to \infty} \lambda \Vol(M_\tau \cap W_\ell) = \ell L(\Sigma) .
$$
Finally, by (\ref{eq:Mtau-Well}) we have 
$$
\lim_{\lambda \to \infty} \lambda \Vol(M_\tau) = \ell L(\Sigma) .
$$
For $n \geq 3$, working in an $h$--orthonormal frame with $e_{n-1}$
proportional to $D_x u_\tau$, so that $D_x u_\tau = u_{\tau,x} e_{n-1}$, with
$|D_x u_\tau|_h^2 = u_{\tau,x}^2$, we have 
$$
g_\tau = \begin{pmatrix} u_\tau^2 h_{\perp} & 0 & 0 \\
                        0 & - u_{\tau,x}^2  + u_\tau^2 h_{//} & - u_{\tau,x}
                        \partial_r u_\tau \\
                       0 & - u_{\tau,x} \partial_r u_\tau & - (\partial_r
                        u_\tau)^2 + 1 
\end{pmatrix} ,
$$
where $h_{\perp}$, $h_{//}$ the restriction of $h$ to $e_{n-1}^{\perp}$ and
to $e_{n-1}$ respectively. 
This gives 
$$
\det g_\tau = \det (u_\tau^2 h_{\perp} ) \left ( (1-(\partial_r u_\tau)^2 )
u_\tau^2 h_{//} - |D_x u_\tau|_h^2 \right )  .
$$
Taking into account $\det h = (\det h_{\perp} )h_{//}$,
and arguing by analogy with the case $n=2$ shows that 
$$
\lim_{\lambda \to \infty} \lambda^{n-1} \Vol(M_\tau) = \ell \Vol(\Sigma) .
$$
Next we consider the distance function on $M_\tau$. 
Let $p,q \in M_\tau \cap W_\ell$ and let $\gamma$ be a curve
connecting $p,q$. We may restrict our consideration to curves such that
$dr(\dot \gamma) \ne 0$. By parametrizing $\gamma = \gamma(s)$ so that 
$dr(\dot \gamma) = 1$, we have
$$
|\dot \gamma |_{g_\tau}^2 = u_\tau^2 |\dot \gamma_x|_h^2 + 1 - |du_\tau (\dot
 \gamma)|^2 = 1 + O(\lambda^{-2}) ,
$$
and hence 
$$
L[\gamma] = |r(p) - r(q)| + O(\lambda^{-2}) .
$$
We state the conclusions of this section as 
\begin{thm}\label{thm:wedge} 
\begin{enumerate}
\item $\lim_{\tau \to -\infty} \lambda^{n-1} \Vol(M_\tau) = \ell
\Vol(\Sigma)$,
\item As $\tau \to - \infty$, the geometry of $M_{\tau}$ converges in the
Gromov sense to the interval of length $\ell$.
\end{enumerate}
\end{thm}
\mnote{From the point of view of Benedetti and Guadagnini, $M_{\tau}$
converges to the R-tree dual to the measured foliation given by $(\Sigma,
\ell)$}
 
\section{CMC hypersurfaces in simplicial flat spacetimes}\label{sec:simplicial}

Let $(M,g)$ be a compact hyperbolic manifold with metric $g$ of sectional
curvature $-1$, of dimension $n \geq 2$, 
and let $\LL = \{ (\Sigma_k, \ell_k), k = 1,\dots,m\}$ 
be a weighted, finite
collection of nonintersecting compact simple totally geodesic hypersurfaces
with weigths $\ell_k \in \Re$,  in $(M,g)$. Further, let $V$ be the
simplicial flat spacetime obtained by performing elementary deformations
w.r.t. the elements $(\Sigma_k, \ell_k)$ of $\LL$. Let $(\Tree,d)$ be the
simplicial $\Re$--tree dual to $\LL$. 

Let $M_\tau$ be the leaves of the global CMC foliation of $V$, and 
let
$\ell_{M_\tau}, s_{\Tree}$ be the marked length spectrum of $M_\tau$ and the
marked measure spectrum of $\Tree$, respectively. 
The conclusion of Theorem \ref{thm:wedge} generalizes immediately to the
situation of simplicial flat spacetimes. 
\begin{thm}\label{thm:simpconv}
With the notation introduced above, the following holds.
\begin{enumerate}
\item 
$\lim_{\tau \to -\infty} \tilde M_\tau = \Tree$, where the limit is
understood in the Gromov sense. Thus, 
the induced geometry on the universal cover $\tilde M_\tau$ converges, as
$\tau \to -\infty$, in the
Gromov sense to $(\Tree, d)$. 
\item 
$\lim_{\tau \to -\infty} \ell_{M_\tau} = s_{\Tree}$
\end{enumerate}
\end{thm}

\section{Dirichlet energy and rescaled Hamiltonian}
The Gauss map $\varphi: M_{\tau}
\to M$ is harmonic from the CMC hypersurface $M_{\tau}$ with its induced
geometry to $M$ with its hyperbolic geometry \cite{andersson:flatcmc}. 
Further, $\varphi$ is the unique
harmonic map $M_{\tau} \to M$ isotopic to the identity. The harmonic map (Dirichlet) energy
of $\varphi$,  
defined by 
$E(M_{\tau} , \varphi) = \int_{M_\tau} |d\varphi|^2\mu_g$,
can be written as 
$$
E(M_{\tau} , \varphi) 
= \int_M |K|^2 d\mu_g = \int_{M_\tau} R \mu_g 
+ \tau^2 \Vol(M,g) , 
$$
In case $n=2$, $\int_M R \mu_g = 4\pi \chi(M)$ by Gauss-Bonnet, which
gives the interesting formula 
$$
E(M_{\tau}, \varphi) = 4 \pi \chi(M_\tau) + \tau^2 \Vol(M_\tau,g) ,
$$
found by Puzio \cite{puzio:gauss}. 

The rescaled Hamiltonian 
$\Ham = |\tau|^n \Vol(M_\tau)$ is the Hamiltonian for gravity in a suitably
chosen gauge, see \cite{fischer:moncrief:sigma}. As shown in
\cite{andersson:flatcmc} it satisfies
\begin{equation}\label{eq:Hamineq:thm}
\Ham \geq n^n \Vol(M,g) .
\end{equation}
Equality in (\ref{eq:Hamineq:thm}) holds if and only if $V$ is the Lorentz cone
over $(M,g)$. 

We will use our results on flat simplicial space--times to 
understand the limiting behavior of the Dirichlet energy and the rescaled Hamiltonian. 
Let $V,M,\LL$ be as in section \ref{sec:simplicial}, and let $M_\tau,
M_\lambda$ be the leaves of the CMC foliation of $V$ and the rescaled leaves, respectively. 
If we let $\lambda = |\tau|/(n-1)$ as above, then 
\begin{equation}\label{eq:Escale}
E(M_{\lambda}, \varphi) = \lambda^{n-2} E(M_{\tau} , \varphi)
\end{equation}
is scale invariant. Let $K_{\lambda}$ be the second fundamental from of
$M_{\lambda}$. Since we know from the above that for a wedge spacetime, 
the height functions 
$u_{\lambda} \to
w_{\lambda}$, we are able to conclude in the simplicial case, from the bounds on the derivatives of
$u_{\lambda}$ that $K_{\lambda} \to h \oplus 0$, on each wedge, 
where $h$ is the metric on
$\Sigma$. We have $|h \oplus 0|^2 = n-1$, which taking into account the fact
that the contribution from the part of $M_\lambda$ off the wedge can be
ignored by the arguments above, gives that the contribution from each wedge
to $E(M_{\lambda} , \varphi)$ behaves
like $(n-1)\Vol(M_{\lambda})$. By (\ref{eq:Escale}), this gives 
$\lim_{\lambda \to \infty} \lambda^{n-3} E(M_{\tau} , \varphi) =
(n-1)\sum_{k=1}^m \ell_k  
\Vol(\Sigma_k)
$. 

The rescaled Hamiltonian $\Ham$ is scale invariant, so we can consider its
behavior on $M_\lambda$. Here we have mean curvature approximate to 1 in the
wedges, while the wedges have length $\lambda \ell_k$, which gives 
$\lim_{\lambda \to \infty} \lambda^{-1} \Ham = \sum_{k=1}^m \ell_k
\Vol(\Sigma_k)$. 
Summarizing, we have   
\begin{thm}\label{thm:vol} \begin{enumerate}
\item 
$\lim_{\lambda \to \infty} \lambda^{n-3} E(M_{\tau} , \varphi) =
(n-1)\sum_{k=1}^m \ell_k  
\Vol(\Sigma_k)
$,
\item
$\lim_{\lambda \to \infty} \lambda^{-1} \Ham = \sum_{k=1}^m \ell_k
\Vol(\Sigma_k)$. 
\end{enumerate} 
Specializing to the $2+1$ dimensional case, we have denoting the length of
$\Sigma$ by $L(\Sigma)$, 
$
\lim_{\lambda \to \infty} \lambda^{-1} E(M_{\tau} , \varphi) =
\sum_{k=1}^m \ell_k 
L(\Sigma_k),$ and similarly for $\Ham$. 
\end{thm}
Let us compare this result to what is known about the time dependence of the
Dirichlet energy in the $2+1$ dimensional case. It has been proved 
\cite[Lemma 4.4]{andersson:etal:2+1grav} that with our present conventions,
\begin{align*}
A(M_{\tau_0})|\tau_0|^2/|\tau| \leq |\tau|A(M_{\tau}) &\leq
|\tau_0|A(M_{\tau_0}), \qquad \text{\rm for } \tau < \tau_0 < 0 \\
A(M_{\tau_0})|\tau_0| \leq |\tau|A(M_{\tau}) &\leq
|\tau_0|^2 A(M_{\tau_0})/|\tau| , \qquad \text{\rm for } \tau_0 < \tau < 0 
\end{align*}
This together with Puzio's result gives 
$$
A(M_{\tau_0})|\tau_0|^2 \leq E(M_\tau, \varphi) - 4 \pi \chi(M_{\tau}) \leq 
A(M_{\tau_0}) ) |\tau| |\tau_0| ,
\quad \text{ for $\tau < \tau_0 < 0$}
$$
which gives the correct leading order behavior in the collapsing direction, 
but which does not identify
the coefficient. 
Similarly in the expanding direction we know that 
$$
\lim_{\lambda \to 0} \lambda^2 A(M_{\tau}) = A(M,g)
$$
the area of the hyperbolic geometry on $M$. Therefore we find that 
$$
\lim_{\tau \to 0} E(M_\tau,\varphi) = 4 \pi \chi(M_{\tau}) + 4 A(M,g) 
= 4 \pi | \chi(M)|.
$$
\bigskip

\noindent{\bf Acknowledgements:} The author is grateful Vince Moncrief,
Ralph Howard, Riccardo Benedetti and Kevin Scannell
for helpful discussions
on various aspects of this material.


\begin{thebibliography}{10}

\bibitem{andersson:flatcmc}
Lars Andersson, \emph{Constant mean curvature foliations of flat space-times},
  Comm. Anal. Geom. \textbf{10} (2002), no.~5, 1125--1150.

\bibitem{andersson:etal:cosmtime}
Lars Andersson, Gregory~J. Galloway, and Ralph Howard, \emph{The cosmological
  time function}, Classical Quantum Gravity \textbf{15} (1998), no.~2,
  309--322.

\bibitem{andersson:etal:maxprin}
\bysame, \emph{A strong maximum principle for weak solutions of quasi-linear
  elliptic equations with applications to {L}orentzian and {R}iemannian
  geometry}, Comm. Pure Appl. Math. \textbf{51} (1998), no.~6, 581--624.

\bibitem{andersson:etal:2+1grav}
Lars Andersson, Vincent Moncrief, and Anthony~J. Tromba, \emph{On the global
  evolution problem in $2+1$ gravity}, J. Geom. Phys. \textbf{23} (1997),
  no.~3-4, 191--205.

\bibitem{benedetti:guadagnini:cosm}
Riccardo Benedetti and Enore Guadagnini, \emph{Cosmological time in
  {$(2+1)$}-gravity}, Nuclear Phys. B \textbf{613} (2001), no.~1-2, 330--352.

\bibitem{cheng:yau:maximal}
Shiu~Yuen Cheng and Shing~Tung Yau, \emph{Maximal space-like hypersurfaces in
  the {L}orentz-{M}inkowski spaces}, Ann. of Math. (2) \textbf{104} (1976),
  no.~3, 407--419.

\bibitem{treibergs:choi}
Hyeong~In Choi and Andrejs Treibergs, \emph{Gauss maps of spacelike constant
  mean curvature hypersurfaces of {M}inkowski space}, J. Differential Geom.
  \textbf{32} (1990), no.~3, 775--817.

\bibitem{fischer:moncrief:sigma}
Arthur~E. Fischer and Vincent Moncrief, \emph{The {E}instein flow, the
  $\sigma$-constant and the geometrization of $3$-manifolds}, Classical Quantum
  Gravity \textbf{16} (1999), no.~11, L79--L87.

\bibitem{gerhardt:CMC}
Claus Gerhardt, \emph{${H}$-surfaces in {L}orentzian manifolds}, Comm. Math.
  Phys. \textbf{89} (1983), no.~4, 523--553.

\bibitem{kapovich:deform}
Michael Kapovich, \emph{Deformations of representations of discrete subgroups
  of ${\rm {s}{o}}(3,1)$}, Math. Ann. \textbf{299} (1994), no.~2, 341--354.

\bibitem{kapovich:millson:deform}
Michael Kapovich and John~J. Millson, \emph{On the deformation theory of
  representations of fundamental groups of compact hyperbolic $3$-manifolds},
  Topology \textbf{35} (1996), no.~4, 1085--1106.

\bibitem{mess:const:curv}
Gerhard Mess, \emph{Lorentz spacetimes of constant curvature}, Tech. Report
  IHES/M/90/28, Institute des Hautes Etudes Scientifiques, 1990.

\bibitem{otal:hyp}
Jean-Pierre Otal, \emph{Le th\'eor\`eme d'hyperbolisation pour les vari\'et\'es
  fibr\'ees de dimension 3}, Ast\'erisque (1996), no.~235, x+159.

\bibitem{puzio:gauss}
Raymond~S. Puzio, \emph{The {G}auss map and $2+1$ gravity}, Classical Quantum
  Gravity \textbf{11} (1994), no.~11, 2667--2675.

\bibitem{scannell:deform}
Kevin~P. Scannell, \emph{Infinitesimal deformations of some {${\rm
  {S}{O}}(3,1)$} lattices}, Pacific J. Math. \textbf{194} (2000), no.~2,
  455--464.

\bibitem{schoen:yau:incompressible}
Richard~M. Schoen and Shing~Tung Yau, \emph{Existence of incompressible minimal
  surfaces and the topology of three-dimensional manifolds with nonnegative
  scalar curvature}, Ann. of Math. (2) \textbf{110} (1979), no.~1, 127--142.

\bibitem{treibergs:cmc}
Andrejs~E. Treibergs, \emph{Entire spacelike hypersurfaces of constant mean
  curvature in {M}inkowski space}, Invent. Math. \textbf{66} (1982), no.~1,
  39--56.

\bibitem{tromba:book}
Anthony~J. Tromba, \emph{Teichm\"uller theory in {R}iemannian geometry},
  Lectures in Mathematics ETH Z\"urich, Birkh\"auser Verlag, Basel, 1992,
  Lecture notes prepared by Jochen Denzler.

\bibitem{zeghib:hyp:geod}
Abdelghani Zeghib, \emph{Laminations et hypersurfaces g\'eod\'esiques des
  vari\'et\'es hyperboliques}, Ann. Sci. \'Ecole Norm. Sup. (4) \textbf{24}
  (1991), no.~2, 171--188.

\end{thebibliography}

\providecommand{\bysame}{\leavevmode\hbox to3em{\hrulefill}\thinspace}
\providecommand{\MR}{\relax\ifhmode\unskip\space\fi MR }
\providecommand{\MRhref}[2]{%
  \href{http://www.ams.org/mathscinet-getitem?mr=#1}{#2}
}
\providecommand{\href}[2]{#2}

\end{document}